\documentclass[11pt]{article}
\usepackage{url}
\usepackage[all]{xy}

\usepackage{amssymb,amsbsy,latexsym,amsmath,bbm,epsfig,psfrag,amsthm,mathrsfs,stmaryrd}
\usepackage[swedish,english]{babel}
\usepackage{graphicx,tikz,esint}
\usepackage[T1]{fontenc}
\usepackage[latin1]{inputenc}
\usepackage{amsfonts}
\usepackage{amsmath}
\usepackage{amssymb}
\usepackage{epsf}
\usepackage{graphics}
\usepackage{graphicx}
\usepackage{latexsym}
\usepackage{psfrag}
\usepackage{relsize}
\usepackage{verbatim}
\usepackage{hyperref}
\usepackage{pgfplots}
\usepackage{float}
\usetikzlibrary{positioning}
\usetikzlibrary{arrows}
\usetikzlibrary{decorations.pathreplacing}

\makeatletter
\DeclareFontFamily{OMX}{MnSymbolE}{}
\DeclareSymbolFont{MnLargeSymbols}{OMX}{MnSymbolE}{m}{n}
\SetSymbolFont{MnLargeSymbols}{bold}{OMX}{MnSymbolE}{b}{n}
\DeclareFontShape{OMX}{MnSymbolE}{m}{n}{
    <-6>  MnSymbolE5
   <6-7>  MnSymbolE6
   <7-8>  MnSymbolE7
   <8-9>  MnSymbolE8
   <9-10> MnSymbolE9
  <10-12> MnSymbolE10
  <12->   MnSymbolE12
}{}
\DeclareFontShape{OMX}{MnSymbolE}{b}{n}{
    <-6>  MnSymbolE-Bold5
   <6-7>  MnSymbolE-Bold6
   <7-8>  MnSymbolE-Bold7
   <8-9>  MnSymbolE-Bold8
   <9-10> MnSymbolE-Bold9
  <10-12> MnSymbolE-Bold10
  <12->   MnSymbolE-Bold12
}{}

\let\llangle\@undefined
\let\rrangle\@undefined
\DeclareMathDelimiter{\llangle}{\mathopen}%
                     {MnLargeSymbols}{'164}{MnLargeSymbols}{'164}
\DeclareMathDelimiter{\rrangle}{\mathclose}%
                     {MnLargeSymbols}{'171}{MnLargeSymbols}{'171}
\makeatother

\theoremstyle{plain}

\theoremstyle{definition}
\newtheorem{Lem}{Lemma}
\numberwithin{Lem}{section}
\newtheorem{Prop}{Proposition}
\numberwithin{Prop}{section}
\newtheorem{Thm}{Theorem}
\numberwithin{Thm}{section}

\numberwithin{Cor}{section}

\numberwithin{Con}{section}

\numberwithin{Prob}{section}

\newtheorem{Def}{Definition}
\numberwithin{Def}{section}

\numberwithin{hyp}{section}

\numberwithin{conj}{section}
\newtheorem{ex}{Example}
\numberwithin{ex}{section}

\theoremstyle{remark}
\newtheorem{rem}{\bf{Remark}}
\numberwithin{rem}{section}

\numberwithin{equation}{section}

\newcommand{\dv}{\partial}
\newcommand{\Om}{\Omega}

\newcommand{\eps}{\varepsilon}
\newcommand{\R}{{\mathbb R}}
\newcommand{\C}{{\mathbb C}}

\newcommand{\Di}{\mathbb{D}}

\newcommand{\A}{\mathcal{A}}

\begin{document}

\vspace{.4cm}

\title{\sffamily Conformal Structure of Autonomous Leray-Lions Equations in the Plane and Linearisation by Hodograph Transform}
\author{}
\date{}
\maketitle

\begin{abstract}
We give sufficient conditions for when an autonomous elliptic Leray-Lions equation in the plane has a conformal structure. This allows the Leray-Lions equation to be linearised in a special form through the hodograph transform. 
\end{abstract}

\subsection*{Keywords}
Elliptic partial differential equations, Beltrami equation, conformal structure

\section{\sffamily Introduction and Motivation}

Consider a general autonomous second order equation in the plane of the form 
\begin{align}\label{eq:LL}
\text{div}\, \A(\nabla u(z))=0, \quad z\in \Om
\end{align}
for some domain $\Om\subset \C$ and some continuous monotone field $\A\in W^{1,2}_{loc}(\R^2,\R^2)\cap C(\R^2,\R^2)$, whose precise assumption we defer to Definition \ref{def:Dmon}. In particular
\begin{align*}
\langle \A(\xi)-\A(\zeta),\xi-\zeta\rangle > 0
\end{align*}
for all $\xi\neq \zeta\in \C$. These are very weak assumptions on $\A$, and it implies that the equation can be a highly degenerate elliptic equation. Assume that $u$ is a weak solution of \eqref{eq:LL} such that $\nabla u\in W^{1,2}_{loc}(\Om,\R^2)$. Then by looking at the complex gradient $f=u_z$ one can show (and we will recall later how in section \ref{sec:Cgrad}) that the complex gradient solves the Beltrami equation 
\begin{align}\label{eq:Cgrad}
f_{\overline{z}}=\overline{\nu(f)}f_z+\nu(f)\overline{f_z}. 
\end{align}
where 
\begin{align*}
\nu(f):=-\frac{\mathbf{A}_{\overline{w}}(f)}{2\mathfrak{R}e[\mathbf{A}_{w}(f)]}
\end{align*}
and 
\begin{align*}
\mathbf{A}(\xi)=\overline{\A(2\overline{\xi})}. 
\end{align*}

The resulting Beltrami equation for the complex gradient is uniformly elliptic if and only if 
\begin{align*}
2\vert \nu(w)\vert\leq k<1. 
\end{align*}
for all $w\in \C$. Interestingly, one can show that the equation for the complex gradient is uniformly elliptic even though the structure field $\A$ \emph{does not satisfy} the standard uniform ellipticity condition 
\begin{align}\label{def:ellip}
\vert \xi\vert^2+\vert \A(\xi)\vert^2\leq \bigg(K+\frac{1}{K}\bigg)\langle \A(\xi),\xi\rangle
\end{align}
for all $\xi\in \C$ and some $K\geq 1$. This happens in particular for the $p$-Laplacian when $p\neq 2$, which shows that this notion of uniform ellipticity is distinct form \eqref{def:ellip}. Assuming that the Beltrami equation is uniformly elliptic on some domain $\Om\in \C$ it follows that $f$ is $K$-quasiregular for some $1\leq K<+\infty$. This means that $f$ belongs to $W^{1,2}_{loc}(\Om,\C)$ and also solves a $\C$-linear Beltrami equation 
\begin{align*}
f_{\overline{z}}(z)=\mu(z)f_z(z),
\end{align*}
but where $\mu$ \emph{depends on $f$ and is different for different solutions of \eqref{eq:Cgrad}}. In particular the Stoilow factorization theorem implies that $f=h(\chi)$ for some holomorphic function $h$ and a quasiconformal map $\chi$. In many instances this is enough if one only wants to deduce interior regularity of solutions as in \cite{IM}.

More precisely, assume that we insert $f=h\circ \chi$ into \eqref{eq:Cgrad}, with $\chi$ being a quasiconformal homeomorphism and $h$ holomorphic, then we get the equation 
\begin{align*}
\chi_{\overline{z}}(z)=\overline{\nu(h(\chi(z)))}\chi_z(z)+\nu(h(\chi(z)))\frac{\overline{h_w(\chi(z))}}{h_w(\chi(z))}\overline{\chi_z(z)}.
\end{align*}

Applying the hodograph transform gives the equation for the inverse map $\eta=\chi^{-1}$
\begin{align*}
\eta_{\overline{w}}(w)=-\overline{\nu(h(w))}\overline{\eta_w(w)}-\nu(h(w))\frac{\overline{h_w(w)}}{h_w(w)}\eta_w(w)
\end{align*}
which is a linear equation of course, but whose coefficients depend on the holomorphic map and \emph{its derivative $h_w$}. In some instances, especially when one wants to consider boundary behaviour and highly degenerate elliptic equations, this may cause difficulties.  

We are however in no way restricted to only considering the complex gradient $u_z$ of solutions \eqref{eq:LL}, but we could also consider complex fields of the form $F(z)=\Phi(u_z(z))$ for some homeomorphism $\Phi\in W^{1,2}_{loc}(\C,\C)$. If we could find a $\Phi$ so that $F$ solves an autonomous $\C$-quasilinear equation, then when linearising using the Stoilow factorization and the hodograph transform we would get a linear equation whose coefficients depend on a holomorphic function but \emph{not on its derivative}. In particular in the works \cite{ADPZ,K} this feature was decisive.

 Moreover, if in addition solutions to \eqref{eq:LL} are $C^1$ but the structure field $\A\in C^1$ does not satisfy 
\begin{align*}
\bigg\Vert \frac{\mathbf{A}_{\overline{w}}(w)}{2\mathfrak{R}e[\mathbf{A}_{w}(w)]}\bigg\Vert_\infty=k<1
\end{align*}
we could consider a connected component $N$ of
\begin{align*}
\bigg\{z\in \C:\bigg\vert \frac{\mathbf{A}_{\overline{w}}(w)}{2\mathfrak{R}e[\mathbf{A}_{w}(w)]}\bigg\vert_\infty<1\bigg\}
\end{align*} 
and solutions $u$ of \eqref{eq:LL} and open subsets $U\subset \Om$ such that $f=u_z:U\to N$.  
We will now consider a number of instructive examples where this is possible.

\begin{ex}[$p$-Laplacian]

The $p$-Laplacian is the equation
\begin{align*}
\text{div}\, \vert \nabla u(z)\vert^{p-2}\nabla u(z)=0
\end{align*}
In this case $\A(\xi)=\vert \xi\vert^{p-2}\xi$ and one can see that the structure field $\A$ is not uniformly elliptic. Yet, the complex gradient $f=u_z$  solves 
\begin{align*}
f_{\overline{z}}(z)=\bigg(\frac{1}{p}-\frac{1}{2}\bigg)\bigg[\frac{\overline{f}}{f}f_z+\frac{f}{\overline{f}}\overline{f_z}\bigg],
\end{align*}
and so 
\begin{align*}
\mu(f)=\bigg(\frac{1}{p}-\frac{1}{2}\bigg)\frac{\overline{f}}{f}, \quad \nu(f)=\bigg(\frac{1}{p}-\frac{1}{2}\bigg)\frac{f}{\overline{f}}.
\end{align*}
which is uniformly elliptic! Let $\Phi_\delta(z)=\vert z\vert^{\delta-1}z$ where $\delta=\sqrt{p-1}$. Then it is shown in \cite[ch. 16]{AIM} that $F=\Phi_\delta(f)$ solves the $\C$-quasilinear Beltrami equation 
\begin{align*}
F_{\overline{z}}=\frac{1-\delta}{1+\delta}\frac{\overline{F}}{F}F_z. 
\end{align*}

Using the Stoilow factorization theorem $F=\phi\circ \chi$ for some holomorphic $\phi$ and a homeomorphic solution $\chi$ of 
\begin{align*}
\chi_{\overline{z}}=\frac{1-\delta}{1+\delta}\frac{\overline{F}}{F}\chi_z. 
\end{align*}
If we let $g(z)=\chi^{-1}(z)$, then the hodograph transform gives the linear equation 
\begin{align*}
g_{\overline{z}}(z)=-\frac{1-\delta}{1+\delta}\frac{\overline{\phi(z)}}{\phi(z)}\overline{g_z(z)}
\end{align*}
\end{ex}

\begin{ex}[Minimal surfaces]
\label{ex:Mini}

In this case the autonomous Leray-Lions equation equals 

\begin{align*}
\text{div}\,\frac{\nabla u(z)}{\sqrt{1+\vert \nabla u(z)\vert^2}}=0
\end{align*}
and $\A(\xi)=\frac{\xi}{\sqrt{1+\vert \nabla \xi\vert^2}}$. 
If one lets 
\begin{align*}
\Phi(z)=\frac{2z}{1+\sqrt{1+4\vert z\vert^2}}
\end{align*}
then it was shown in \cite{K} (see also \cite[Lemma 5.1., p.169]{KS}) that $F=\Phi(u_z)$ solves the $\C$-quasilinear Beltrami equation 
\begin{align*}
F_{\overline{z}}=\overline{F}^2F_z. 
\end{align*}
Using the Stoilow factorisation, $F=\phi\circ \chi$ where $\chi$ is a homeomorphic solution to 
\begin{align*}
\chi_{\overline{z}}=\overline{F}^2\chi_z,
\end{align*}
the hodograph transform then gives that $g=\chi^{-1}$ solves the linear equation
\begin{align*}
g_{\overline{z}}(z)=-\overline{\phi(z)}^2\overline{g_z(z)}.
\end{align*}
\end{ex}

\begin{ex}[Dimer models]
When studying the asymptotic behaviour of random height functions in dimer models, e.g. \cite{KOS,CKP,ADPZ} one is lead to the study of the Euler-Lagrange equation 
\begin{align*}
\text{div}\, \nabla \sigma(\nabla u(z))=0
\end{align*}
where $\sigma$ is an in general inexplicit convex function that solves a boundary value problem for the Monge-AmpÞre equation, see \cite{ADPZ}. By considering the \emph{Lewy transform}
\begin{align*}
L_\sigma(z)=z+\nabla \sigma(z),
\end{align*} 
it is shown in \cite{ADPZ} that the complex valued field 
\begin{align*}
F(z)=\overline{L_\sigma(\nabla u(z))}=\overline{L_\sigma( 2\overline{_zu(z)})}
\end{align*}
solves a $\C$-quasilinear equation of the form 
\begin{align*}
F_{\overline{z}}(z)=\mathcal{H}'(F(z))F_z(z),
\end{align*}
where $\mathcal{H}'$ is a proper holomorphic map. The Stoilow factorisation and the hodograph transform then shows that for $F=\phi\circ \chi$, with $\phi$ holomorphic and $\chi$ a homeomorphic solution of 
\begin{align*}
\chi_{\overline{z}}=\mathcal{H}'(F)\chi_z
\end{align*}
the inverse $g=\chi^{-1}$ solves the linear Beltrami equation 
\begin{align*}
g_{\overline{z}}(z)=-\mathcal{H}'(\phi(z))\overline{g_z(z)}. 
\end{align*}
\end{ex}

Given these examples one can ask if it is always possible to find a homeomorphism $\Phi$ such that the complex field $F=\Phi(u_z)$ solves an elliptic $\C$-quasilinear Beltrami equation if $u$ is a solution of \eqref{eq:LL}?

Assuming for now the existence of such a $\Phi$ it follows that $F$ would solve a $\C$-quasilinear equation of the form 
\begin{align}\label{eq:CB1}
F_{\overline{z}}(z)=\gamma(F(z))F_{z}(z).
\end{align}
If we assume that $u$ solves \eqref{eq:LL} on an open set $U\subset \C$ we also need to assume here, to get a viable theory, that $\vert \gamma(F(z))\vert\leq k_V<1$ on any relatively compact set $V\subset U$. If this holds we call $\gamma$ \emph{the conformal Beltrami coefficient} associated the structure field $\A$. 

We may now apply the Stoilow factorization which says that every solution $F$ of \eqref{eq:CB1} is of the form 
\begin{align*}
F= h\circ \chi
\end{align*}
where $h$ is a holomorphic function and $\chi$ is a homeomorphic solution to 
\begin{align*}
\chi_{\overline{z}}(z)=\gamma(h(\chi(z)))\chi_{z}(z)
\end{align*}

If we let $g=\chi^{-1}$ the hodograph transformation yields (see \cite[ch. 16.3]{AIM}) that $g$ solves the \emph{anti-$\C$-linear Beltrami equation}
\begin{align}\label{eq:LinB}
g_{\overline{z}}(z)=-\gamma(h(z))\overline{g_{z}(z)}
\end{align}

In particular, if the regularity of $\gamma$ is known so is $\gamma\circ h$ and \emph{does not depend} on $g$ itself which is major advantage. Moreover, the equation degenerates precisely when $\vert \gamma \circ h(z)\vert=1$.


\section{\sffamily Main results}

\begin{Def}\label{def:EtaGamma}
Let $N\subset \C$ be open and let $\A: N\to \C$ be a locally $\delta$-monotone field according to Definition \ref{def:Dmon}. Let $\mathbf{A}(z)=\overline{\A(2\overline{z})}$ and define 
\begin{align}
\nu(z):=-\frac{\mathbf{A}_{\overline{z}}(z)}{2\mathfrak{R}e[\mathbf{A}_{z}(z)]}.
\end{align}
Define 
\begin{align}\label{def:EtaGamma}
\eta(z)&=\left\{
    \begin{array}{ll}
    0, & \text{if $\nu(z)=0$},\\
   \displaystyle \frac{-1+\sqrt{1-4\vert \nu(z)\vert^2}}{2\overline{\nu(z)}}& \text{otherwise},
   \end{array} \right. 
\quad \text{and} \quad \gamma(z)=-\overline{\eta(z)}. 
\end{align}
\end{Def}

We now come to the main results of this paper:

\begin{Thm}\label{thm:Main0}
Let $\A: \C\to \C$ be a $\delta$-monotone field (see Definition \ref{def:Delta}) and let $\nu$ and $\gamma$ be given by Definition \ref{def:EtaGamma}. Let $\Om \subset \C$ be a bounded domain and consider all weak solutions $u\in \mathbb{W}^{1,\Upsilon}(\Om)$ of the autonomous Leray-Lions equations  
\begin{align*}
\text{div}\,\A(\nabla u(z))=0, \quad z\in \Om, 
\end{align*} 
where $\mathbb{W}^{1,\Upsilon}(\Om,\C)$ is the homogenous Orlicz-Sobolev space associated to $\A$ as in \cite[ch. 16.4.1]{AIM}.
Then there exists a quasiconformal map $\Phi: \C\to \C$ that solves the uniformly elliptic linear Beltrami equation 
\begin{align*}
\Phi_{\overline{z}}(z)=\eta(z)\Phi_z(z)
\end{align*}
where $\eta$ is given by \eqref{def:EtaGamma} such that the associated field 
\begin{align*}
F(z)=\Phi(u_z)
\end{align*}
solves the uniformly elliptic $\C$-quasilinear equation 
\begin{align}\label{eq:Main}
F_{\overline{z}}(z)=\gamma(\Phi^{-1}(F(z))F_z(z)
\end{align}
for a.e. $z\in \Om$. In particular, $\Phi$ is \emph{independent of $u$} and only depends on $\nu$.  
\end{Thm}

\begin{rem}
By \cite[Theorem 16.4.5]{AIM} any weak solution in $\mathbb{W}^{1,\Upsilon}(\Om)$ of the Leray-Lions equation \eqref{eq:LL} is $C^{1,\alpha}$, in particular the complex gradient $u_z$ is continuous. In addition $\mathbb{W}^{1,\Upsilon}(\Om)\subset W^{1,p}_{loc}(\Om)$ for $p=\frac{2}{1+\sqrt{1-\delta^2}}>1$. 
\end{rem}

Before we state the next theorem we recall the definition of a Koebe domain. 

\begin{Def}
A Koebe domain, also called \emph{circle domain}, is a planar domain such that each connected component of its boundary is either a circle or a point. 
\end{Def}

\begin{Thm}\label{thm:Main}
Let $N\subset \C$ be a finitely connected domain conformally equivalent to a bounded Koebe domain and and let $\A: N\to \C$ be a locally $\delta$-monotone field according to Definition \ref{def:Dmon}.  Furthermore, let $\nu$ and $\gamma$ be given as in Theorem \ref{thm:Main0}.
Let $\Om \subset \C$ be a domain and consider all weak solutions $u$ of the autonomous Leray-Lions equations such that in addition $u_z\in C^1(\Om)\cap W^{1,2}_{loc}(\Om)$ and  $u_z:\Om \to N$. Then there exists a homeomorphism (of finite distortion) $\Phi: N\to N$ independent of $u$ that solves the linear Beltrami equation 
\begin{align*}
\Phi_{\overline{z}}(z)=\eta(z)\Phi_z(z)
\end{align*}
where $\eta$ is given by \eqref{def:EtaGamma} such that the associated field 
\begin{align*}
F(z)=\Phi(u_z)
\end{align*}
solves the $\C$-quasilinear equation 
\begin{align}\label{eq:Main}
F_{\overline{z}}(z)=\gamma(\Phi^{-1}(F(z))F_z(z)
\end{align}
for a.e. $z\in \Om$. 
\end{Thm}

The proof of Theorem \ref{thm:Main0} and Theorem \ref{thm:Main} will be given in several steps in Section \ref{sec:Proof}.

\begin{rem}
The same result as Theorem \ref{thm:Main} was proven in the special case in \cite[Prop 2.1, Thm. 2.2 ]{KP} when $\A=\nabla \sigma$ and $\sigma$ is a smooth strictly convex function on the interior of $N$ using a different proof. However, they make a different choice of homeomorphism $\Phi$, with their $\Phi$ being \emph{orientation reversing}. In section \ref{sec:CounterEx} we give a counter example, showing that one cannot always chose $\Phi$ to be orientation reversing when $\sigma$ is not real analytic. 
\end{rem}

\begin{Def}[Conformal structure]
Let $\A$ and $\Phi$ by as in either Theorem \ref{thm:Main0} or Theorem \ref{thm:Main}. The conformal Beltrami coefficient $\mu$ associated to the Leray-Lions equation \eqref{eq:LL} is defined according to 
\begin{align}\label{def:ConfB}
\mu=\gamma\circ \Phi^{-1}. 
\end{align}
The \emph{conformal} structure of the Leray-Lions equation is given by 
\begin{align}\label{def:ConformalS}
\mu(F(z))=\gamma(u_z). 
\end{align}
\end{Def}

In fact the proofs of Theorem \ref{thm:Main0} and Theorem \ref{thm:Main} show that the conformal Beltrami coefficient $\mu$ is unique up to post-composition by conformal maps. Moreover it is independent of $u$. Note that on the other the conformal structure \emph{depends} on $u$.

\begin{Prop}\label{prop:Main}
Let $\Om\subset N$ be a finitely connected bounded domain. Let $N,\eta,\Phi$ and $F$ be given as in Theorem \ref{thm:Main}. Let $\mathcal{D}$ be a bounded Koebe domain conformally equivalent to $\Om$. Then every non-constant continuous bounded solution $F\in W^{1,2}_{loc}$ of 
\begin{align*}
F_{\overline{z}}(z)=\mu(F(z))F_z(z)
\end{align*}
factorises according to $F=\varphi\circ g^{-1}$ where $g:\mathcal{D}\to \Om$ is a homeomorphic solution of the \emph{linear} equation 
\begin{align}\label{eq:AntiCLin}
g_{\overline{z}}=-\mu(\varphi(z))\overline{g_z(z)}
\end{align}
and $\varphi: \mathcal{D}\to N$ is a holomorphic function. In particular the Beltrami coefficient of in the linear equation for $g$ depends only on a holomorphic function $h$ but \emph{not on its derivative}. 
\end{Prop}

\begin{proof}
The proof is a direct generalisation of the first part of the proof of \cite[Thm. 4.1]{ADPZ}. See also the diagram below. 
\end{proof}

\[
\xymatrix{
  \mathcal{D} \ar[dd]_{g}  \ar[rr]^\varphi & &\,\,\,N \circlearrowleft \Phi\\ &&\\
 \Om \ar[uurr]_{F=\Phi(u_z)} & &}
\]

\begin{Thm}\label{Thm:Main2}
Let $N,\eta,\Phi,\mathcal{D},\mu$ and $\varphi$ be given as in Theorem \ref{thm:Main} and Proposition \ref{prop:Main}. Assume that $\phi:=\mu\circ \varphi\in W^{1,p}(\mathcal{D},N)$ for some $p>2$. Then every solution of \eqref{eq:AntiCLin}, is of the form 
\begin{align*}
g(z)=\frac{1}{1-\vert \phi(z)\vert^2}(h(z)-\phi(z)\overline{h(z)})
\end{align*}
where $h$ is a generalised analytic function in the sense of Vekua and solves the equation 
\begin{align}\label{eq:Vekua}
h_{\overline{z}}(z)=\frac{\phi_{\overline{z}}(z)}{1-\vert \phi(z)\vert^2}(\overline{h(z)}-\overline{\phi(z)}h(z))
\end{align}

In addition every solution of \eqref{eq:Vekua} is of the form 
\begin{align*}
h(z)=\psi(z)e^{\omega(z)}
\end{align*}
where $\psi$ is holomorphic and $\omega$ is given by \eqref{eq:Omega}. 
\end{Thm}

The proof is given in subsection \ref{sec:AntiC}.


\section{\sffamily Complex gradient method}\label{sec:Cgrad}

In this section, for the convenience of the reader, we will survey the complex gradient method and how it is applied to autonomous Leray-Lions equations. This will provide the necessary background for section \ref{sec:Proof}.  More details about the complex gradient method can be found in \cite[ch. 16]{AIM}. We begin by recalling the concept of $\delta$-monotonicity from \cite{Kov}. 
\begin{Def}\label{def:Delta}
A mapping $\A:\C\to \C$ is \emph{$\delta$-monotone} if there exists a
that there exits a $0<\delta\leq 1$ such that 
\begin{align*}
\langle \A(\xi)-\A(\zeta),\xi-\zeta\rangle \geq \delta \vert \A(\xi)-\A(\zeta)\vert \vert \xi-\zeta \vert. 
\end{align*}
for all $\zeta,\xi \in \C$. 
\end{Def}

We recall \cite[Theorem 3.11.6]{AIM}.

\begin{Thm}
Let $0<\delta \leq 1$. A mapping $\A\in W^{1,2}_{loc}(\C)$ is $\delta$-monotone if and only if 
\begin{align}\label{eq:dm}
\vert \A_{\overline{z}}(z)\vert+\delta\vert \mathfrak{I}m[\A_z(z)]\vert\leq \sqrt{1-\delta^2}\mathfrak{R}e[\A_z(z)]
\end{align}
for a.e. $z$. In particular, $\A$ is $K$-quasiconformal where 
\begin{align*}
K=\frac{1+\sqrt{1-\delta^2}}{1-\sqrt{1-\delta^2}},
\end{align*}
and where the bound on distortion is sharp.
\end{Thm}

One can use  \cite[Theorem 3.11.6]{AIM} to define a $\delta$-monotone maps on domains as follows.

\begin{Def}\label{def:Dmon}
Let $0<\delta \leq 1$ and assume $N\subset \C$ is a domain. A mapping $\A\in W^{1,2}_{loc}(N,\C)$ is $\delta$-monotone if \eqref{eq:dm} holds for a.e. $z\in N$ and $\A$ is a homeomorphism. A mapping $\A\in W^{1,2}_{loc}(N,\C)$ is \emph{locally $\delta$-monotone} if for every $U\Subset N$ there exists a $\delta=\delta(U)$ such that $\A\vert_U$ is $\delta$-monotone. 
\end{Def}

\begin{rem}
It follows from the proof of \cite[Theorem 3.11.6]{AIM} that any solution of \eqref{eq:dm} is locally injective. If $N$ is convex it then follows that $\A$ is automatically a homeomorphism.
\end{rem}

We now consider a general autonomous second order elliptic equation of the form 
\begin{align*}
\text{div}\, \A(\nabla u(x))=0
\end{align*}
where $\A$ is $\delta$-monotone. We will follow the exposition in \cite[ch. 16.4.3, p. 445-447]{AIM}.

Define the new structure field 
\begin{align*}
\mathbf{A}(\xi)=\overline{\A(2\overline{\xi})}. 
\end{align*}

Then $\mathbf{A}$ is monotone as well. Set $f=u_z$. Then $\A(\nabla u(z))=\A(2\overline{u_z})=\overline{\mathbf{A}(f)}$. Thus 
\begin{align*}
\text{div}\, \A(\nabla u)=\text{div}\, \overline{\mathbf{A}(f)}
\end{align*}
Moreover, $0=\text{curl}\,\nabla u(z)=2\text{curl}\,\overline{u_z}=2\text{curl}\,\overline{f}$. Thus the equation \eqref{eq:LL} becomes equivalent to 

\begin{align*}
\text{div}\, \overline{\mathbf{A}(f)}=0
\end{align*}

and then the above together with $\text{curl}\,\nabla u(z)=0$ becomes equivalent to the system
\begin{align}\label{eq:S1}
 \left\{
    \begin{array}{r}
     \text{div}\, \overline{\mathbf{A}(f(z))}=0,\\
     \text{curl}\,\overline{f(z)}=0. 
    \end{array} \right.
\end{align}

We now recall that for any vector field $v$, 
\begin{align*}
\text{div}\,v(z)&=0\quad \Longleftrightarrow\quad \mathfrak{R}e[\dv_zv(z)]=0\\
\text{curl}\,v(z)&=0\quad \Longleftrightarrow\quad \mathfrak{I}m[\dv_zv(z)]=0. 
\end{align*}

Since $\dv_{z}\overline{ f(z))}=\overline{ \dv_{\overline{z}}f(z)}$ the system \eqref{eq:S1} is equivalent to

\begin{align}\label{eq:S2}
 \left\{
    \begin{array}{r}
    \mathfrak{R}e[ \overline{\dv_{\overline{z}}\mathbf{A}(f(z))}]=0,\\
     \mathfrak{I}m[\overline{\dv_{\overline{z}}f(z)}]=0. 
    \end{array} \right.
\quad \Longleftrightarrow\quad
 \left\{
    \begin{array}{r}
    \mathfrak{R}e[ \dv_{\overline{z}}\mathbf{A}(f(z))]=0,\\
     \mathfrak{I}m[\dv_{\overline{z}}f(z)]=0. 
    \end{array} \right.
\end{align}

By \cite[Theorem 16.4.5]{AIM}, both $f$ and $\mathbf{A}(f)$ belong to $W^{1,2}_{loc}$ and in addition are quasiregular on relatively compact subset of $\Om$. In particular the chain rule applies in the pointwise sense and we get 

\begin{align*}
\dv_{\overline{z}}\mathbf{A}(f)=\mathbf{A}_w(f)f_{\overline{z}}+\mathbf{A}_{\overline{w}}(f)\overline{f_{z}}.
\end{align*}

and hence $ \mathfrak{R}e[ \dv_{\overline{z}}\mathbf{A}(f(z))]=0$ is equivalent to 
\begin{align}\label{eq:S3}
\mathbf{A}_w(f)f_{\overline{z}}+\mathbf{A}_{\overline{w}}(f)\overline{f_{z}}+\overline{\mathbf{A}_w(f)f_{\overline{z}}}+\overline{\mathbf{A}_{\overline{w}}(f)\overline{f_{z}}}=0
\end{align}

Moreover $ \mathfrak{I}m[\dv_{\overline{z}}f(z)]=0$ implies $\overline{\dv_{\overline{z}}f(z)}=\dv_{\overline{z}}f(z)$. Inserting this into \eqref{eq:S3} gives

\begin{align*}
\mathbf{A}_w(f)f_{\overline{z}}+\mathbf{A}_{\overline{w}}(f)\overline{f_{z}}+\overline{\mathbf{A}_w(f)}f_{\overline{z}}+\overline{\mathbf{A}_{\overline{w}}(f)}f_{z}=0
\end{align*}
If $\mathbf{A}_{w}(f)+\overline{\mathbf{A}_{w}(f)}\neq 0$ or equivalently $\mathfrak{R}e[\mathbf{A}_w(f)]\neq 0$ we can solve for $f_{\overline{z}}$ giving 
\begin{align*}
f_{\overline{z}}&=-\frac{\overline{\mathbf{A}_{\overline{w}}(f)}}{\mathbf{A}_{w}(f)+\overline{\mathbf{A}_{w}(f)}}f_z-\frac{\mathbf{A}_{\overline{w}}(f)}{\mathbf{A}_{w}(f)+\overline{\mathbf{A}_{w}(f)}}\overline{f_z}\\
&=-\frac{\overline{\mathbf{A}_{\overline{w}}(f)}}{2\mathfrak{R}e[\mathbf{A}_{w}(f)]}f_z-\frac{\mathbf{A}_{\overline{w}}(f)}{2\mathfrak{R}e[\mathbf{A}_{w}(f)]}\overline{f_z}.
\end{align*}
If on the other hand $\mathfrak{R}e[\mathbf{A}_w(f(z))]=0$, one can argue as follows. The set $\{w\in N: \mathfrak{R}e[\mathbf{A}_w(w)]=0\}$ is a null set by \eqref{eq:dm} or else $\mathbf{A}$ is constant which is a contradiction. Let $\mathcal{Z}_{\mathbf{A}}:=\{w\in N: \mathfrak{R}e[\mathbf{A}_w(w)]=0\}$. Since $f=u_z$ is quasiregular, it follows by Stoilow factorization and \cite[Corollary 3.7.6]{AIM} that $f$ satisfies Lusin condition $\mathcal{N}^{-1}$. Consequently, $\vert f^{-1}(\mathcal{Z}_{\mathbf{A}})\vert=0$, and so $\mathfrak{R}e[\mathbf{A}_w(f(z))]\neq 0$ for a.e. $z$. In particular, if $\mathcal{Z}_{\mathbf{A}}$ is a finite set, then can we use \cite[Lemma 7.7,p. 152]{GT} which implies that any $f\in W^{1,1}$, on the set where $f$ is constant we have $f_z=f_{\overline{z}}=0$ a.e.. Thus on that set $f$ automatically solves any $\R$-quasilinear Beltrami equation of the form 
\begin{align}\label{eq:S5}
f_{\overline{z}}=\overline{\nu(f)}f_z+\nu(f)\overline{f_z}\quad \text{for a.e. $z$}
\end{align}
and we are free to define $\nu(f)$ in whatever way we want provided $\vert \nu(f)\vert<1/2$. We may take $\nu=\frac{1}{4}$ for example. Otherwise, whenever $\nu$ is well-defined we let 
\begin{align*}
\nu(f):=-\frac{\mathbf{A}_{\overline{w}}(f)}{\mathbf{A}_{w}(f)+\overline{\mathbf{A}_{w}(f)}}=-\frac{\mathbf{A}_{\overline{w}}(f)}{2\mathfrak{R}e[\mathbf{A}_{w}(f)]}. 
\end{align*}
On the other hand at those points $w$ for which $\nu(w)$ is not defined, which happens in particular for the $p$-Laplace equation at $w=0$, we can again argue as before to see that if $B$ is the set where $\nu$ is not defined, then $\vert f^{-1}(B)\vert=0$, and so $\nu(f(z))$ is well-defined for a.e. $z$ and we may let $\nu(z)=1/4$ in those cases.  Thus we see that the complex gradient $f=u_z$ solves the $\R$-quasilinear Beltrami equation
\begin{align*}
f_{\overline{z}}=\overline{\nu(f)}f_z+\nu(f)\overline{f_z}. 
\end{align*}

The equation is uniformly elliptic if and only if 
\begin{align*}
2\vert \nu(w)\vert\leq k<1. 
\end{align*}
for all $w\in \C$. This holds if $\A$ is $\delta$-monotone on $N$. Otherwise, the equations is uniformly elliptic on relatively compact subsets if $\A$ is locally $\delta$-monotone.   


\section{\sffamily Reduction to $\C$-quasilinear equation and linearisation for the complex gradient equation}\label{sec:Proof}


\subsection{\sffamily Reduction to $\C$-quasilinear equation and linearisation for the complex gradient equation}

Let $N\subset \C$ be an open set and let $\mu,\nu\in C(N,\C)$ and assume that 
\begin{align*}
\vert \mu(w)\vert+\vert \nu(w)\vert<1, \quad w\in N.
\end{align*}

Consider an autonomous $\R$-linear equation 
\begin{align}\label{eq:1}
f_{\overline{z}}=\mu(f) f_z+\nu(f) \overline{f_z},
\end{align}
on a domain $\Om \subset \C$ and assume that all solutions satisfy $f\in C(\Om,\C)\cap W^{1,2}_{loc}(\Om,\C)$ are such that $f:\Om \to N$. Set 
\begin{align*}
F=\Phi(f),
\end{align*}
where $\Phi:N\to N$ is a homeomorphism in $W^{1,2}_{loc}$. Is it possible to choose $\Phi$ such that the new field $F$ solves a $\C$-quasilinear equation? In particular is this possible when $f$ solves the complex gradient equation \eqref{eq:Cgrad}?

\begin{Lem}\label{lem:Chain}
Let $\Phi\in N\to N$ be a homeomorphism in $W^{1,2}_{loc}$ such that for every $U\Subset N$, $\Phi\in R(U,\C)$, where $R(U,\C)$ is the \emph{Royden algebra} of $U$, equal to
\begin{align*}
R(U,\C)=C(U,\C)\cap L^\infty(U,\C)\cap \mathbb{W}^{1,2}(U,\C),
\end{align*}
and where $ \mathbb{W}^{1,2}(U,\C)=\{v\in L^1_{loc}(U,\C):Dv\in L^2\}$ is the homogeneous Sobolev space. 
Let $f$ be a solution of \eqref{eq:1} on $\Om$. Then $F=\Phi(f)$ solves the equation 
\begin{align*}
&(\vert\Phi_w(f)+\Phi_{\overline{w}}(f)\overline{\nu} \vert^2-\vert \Phi_{\overline{w}}(f)\overline{\mu}\vert^2)F_{\overline{z}}\\&=\mu\Big[\vert \Phi_w(f)\vert^2-\vert\Phi_{\overline{w}}(f)\vert^2\Big]F_z+\Big[ \Phi_w(f)^2\nu +\Phi_w(f)\Phi_{\overline{w}}(f)(\vert \nu\vert^2-\vert \mu \vert^2+1)+\Phi_{\overline{w}}(f)^2\overline{\nu}\Big]\overline{F_z}
\end{align*}
\end{Lem}

\begin{proof}
By the assumption on $\Phi$ the chain rule holds, and implies that 
\begin{align}
F_{z}&=\Phi_w(f)f_z+\Phi_{\overline{w}}(f)\overline{f_{\overline{z}}}\\
F_{\overline{z}}&=\Phi_w(f)f_{\overline{z}}+\Phi_{\overline{w}}(f)\overline{f_z}. 
\end{align}
Using \eqref{eq:1} we get 
\begin{align}
F_{z}&=\Phi_w(f)f_z+\Phi_{\overline{w}}(f)(\overline{\mu f_z+\nu \overline{f_z}})\nonumber \\
\label{eq:3}&=(\Phi_w(f)+\Phi_{\overline{w}}(f)\overline{\nu})f_z+\Phi_{\overline{w}}(f)\overline{\mu} \overline{f_z}\\
F_{\overline{z}}&=\Phi_w(f)(\mu f_z+\nu \overline{f_z})+\Phi_{\overline{w}}(f)\overline{f_z}\nonumber\\
\label{eq:4}&=\Phi_w(f)\mu f_z+(\Phi_w(f)\nu+\Phi_{\overline{w}}(f))\overline{f_z}
\end{align}

If we introduce the linear maps $Lw=((\Phi_w(f)+\Phi_{\overline{w}}(f)\overline{\nu}))w+\Phi_{\overline{w}}(f)\overline{\mu}\,\overline{w}$ and $Mw=\Phi_w(f)\mu w+(\Phi_w(f)\nu+\Phi_{\overline{w}}(f))\overline{w}$ we can write the equations as 
\begin{align*}
F_{z}=Lf_z,\quad F_{\overline{z}}=M(f_z)
\end{align*}
and so 
\begin{align*}
F_{\overline{z}}=M\circ L^{-1}(F_z). 
\end{align*}

Using that for a general linear invertible map $Tw=\alpha w+\beta \overline{w}$ the inverse is given by $T^{-1}w=\frac{1}{\vert \alpha\vert^2-\vert \beta\vert^2}(\overline{\alpha}w-\beta \overline{w})$ we get

\begin{align}\label{eq:5}
f_z=\frac{1}{\vert\Phi_w(f)+\Phi_{\overline{w}}(f)\overline{\nu} \vert^2-\vert \Phi_{\overline{w}}(f)\overline{\mu}\vert^2}(\overline{(\Phi_w(f)+\Phi_{\overline{w}}(f)\overline{\nu})}F_z-\Phi_{\overline{w}}(f)\overline{\mu}\overline{F_z})
\end{align}

Inserting this into \eqref{eq:4} gives 
\begin{align*}
&(\vert\Phi_w(f)+\Phi_{\overline{w}}(f)\overline{\nu} \vert^2-\vert \Phi_{\overline{w}}(f)\overline{\mu}\vert^2)F_{\overline{z}}\\
&=\Phi_w(f)\mu (\overline{(\Phi_w(f)+\Phi_{\overline{w}}(f)\overline{\nu})}F_z-\Phi_{\overline{w}}(f)\overline{\mu}\overline{F_z})+(\Phi_w(f)\nu+\Phi_{\overline{w}}(f))((\Phi_w(f)+\Phi_{\overline{w}}(f)\overline{\nu})\overline{F_z}-\overline{\Phi_{\overline{w}}(f)}\mu F_z)\\
&=\mu\Big[\Phi_w(f)\overline{\Phi_w(f)}+\Phi_w(f)\overline{\Phi_{\overline{w}}(f)}\nu-\Phi_w(f)\overline{\Phi_{\overline{w}}(f)}\nu-\Phi_{\overline{w}}(f)\overline{\Phi_{\overline{w}}(f)}\Big]F_z\\
&+\Big[-\vert \mu\vert^2\Phi_w(f)\Phi_{\overline{w}}(f) +(\Phi_w(f)\nu+\Phi_{\overline{w}}(f))(\Phi_w(f)+\Phi_{\overline{w}}(f)\overline{\nu})\Big]\overline{F_z}\\
&=\mu\Big[\vert \Phi_w(f)\vert^2-\vert\Phi_{\overline{w}}(f)\vert^2\Big]F_z\\
&+\Big[ \Phi_w(f)^2\nu +\Phi_w(f)\Phi_{\overline{w}}(f)(\vert \nu\vert^2-\vert \mu \vert^2+1)+\Phi_{\overline{w}}(f)^2\overline{\nu}\Big]\overline{F_z}\\
&=\mu\Big[\vert \Phi_w(f)\vert^2-\vert\Phi_{\overline{w}}(f)\vert^2\Big]F_z\\
&+\Big[ \Phi_w(f)^2\nu +\Phi_w(f)\Phi_{\overline{w}}(f)(\vert \nu\vert^2-\vert \mu \vert^2+1)+\Phi_{\overline{w}}(f)^2\overline{\nu}\Big]\overline{F_z}.
\end{align*}

\end{proof}

Thus $F=\Phi\circ f$ solves a $\C$-quasilinear if and only if 
\begin{align}\label{eq:RedEq}
\Phi_w(f)^2\nu +\Phi_w(f)\Phi_{\overline{w}}(f)(\vert \nu\vert^2-\vert \mu \vert^2+1)+\Phi_{\overline{w}}(f)^2\overline{\nu}=0
\end{align}

provided 

\begin{align}\label{eq:NZ}
\vert\Phi_w(f)+\Phi_{\overline{w}}(f)\overline{\nu} \vert^2-\vert \Phi_{\overline{w}}(f)\overline{\mu}\vert^2\neq 0
\end{align}

We now make the ansatz that a homeomorphic solution of \eqref{eq:NZ} $\Phi$, should it exist, solves a $\C$-linear Beltrami equation
\begin{align}\label{eq:Ansatz}
\Phi_{\overline{z}}(z)=\eta(z)\Phi_z(z)
\end{align}
where the coefficient $\eta$ is to be determined from the equation \eqref{eq:NZ}. This ansatz does not infer any loss of generality since any homeomorphism $\Phi\in W_{loc}^{1,1}(\Om)$ satisfy either 
\begin{align*}
J(z,\Phi)\geq 0 \text{ for a.e. $z\in \Om$, or }\,\,\, J(z,\Phi)\leq 0 \text{ for a.e. $z\in \Om$},
\end{align*}
by \cite[Theorem 3.3.4]{AIM}, where $J(z,\Phi)=\det[D\Phi(z)]=\vert \Phi_z(z)\vert^2-\vert \Phi_{\overline{z}}(z)\vert^2$. If we define 
\begin{align*}
\eta(z):=\frac{ \Phi_{\overline{z}}(z)}{\Phi_z(z)}
\end{align*}
if $\Phi_z(z)\neq 0$ and if $\Phi_z(z)=0$ we define $\eta(z)=0$ if $J(z,\Phi)\geq 0$ a.e. and $\eta(z)=\infty\in \widehat{\C}$ if $J(z,\Phi)\leq 0$ a.e.. Then any homeomorphism solves \eqref{eq:Ansatz}. In particular either $\vert \eta(z)\vert\leq1$ a.e. or $\vert \eta(z)\vert\geq1$.

Inserting this into \eqref{eq:RedEq} gives 
\begin{align}\label{eq:Quad}
&\nu(z)\Phi_z(z)^2 +\Phi_z(z)\Phi_{\overline{z}}(z)(\vert \nu(z)\vert^2-\vert \mu(z) \vert^2+1)+\overline{\nu(z)}\Phi_{\overline{z}}(z)^2\\&=\Bigg[\nu(z) +(\vert \nu(z)\vert^2-\vert \mu(z) \vert^2+1)\eta(z)+\overline{\nu(z)}\eta(z)^2\Bigg]\Phi_z(z)^2=0. 
\end{align}

We impose the condition that
\begin{align}\label{eq:Quad2}
\nu(z) +(\vert \nu(z)\vert^2-\vert \mu(z) \vert^2+1)\eta(z)+\overline{\nu(z)}\eta(z)^2=0. 
\end{align}

which is a quadratic equation provided $\nu(z)\neq 0$.

\begin{Lem}\label{lem:CB1}
Assume $\Phi$ solves the Beltrami equation \eqref{eq:Ansatz} and \eqref{eq:Quad2} holds. Let $\Psi=\Phi^{-1}$. Then $F=\Phi(f)$ solves the $\C$-quasilinear equation 
\begin{align*}
F_{\overline{z}}=\frac{(1-\vert \eta(\Psi(F))\vert^2)\mu(\Psi(F))}{\vert1+\eta(\Psi(F))\overline{\nu(\Psi(F))} \vert^2-\vert \eta(\Psi(F))\overline{\mu(\Psi(F))}\vert^2}F_z.
\end{align*}
\end{Lem}

\begin{proof}
\begin{align*}
F_{\overline{z}}&=\frac{\vert \Phi_w(f)\vert^2-\vert\Phi_{\overline{w}}(f)\vert^2}{\vert\Phi_w(f)+\Phi_{\overline{w}}(f)\overline{\nu(f)} \vert^2-\vert \Phi_{\overline{w}}(f)\overline{\mu(f)}\vert^2}\mu(f)F_z\\
&=\frac{\vert \Phi_w(f)\vert^2-\vert \eta(f)\vert^2\vert\Phi_{w}(f)\vert^2}{\vert\Phi_w(f)+\eta(f)\Phi_{w}(f)\overline{\nu(f)} \vert^2-\vert \eta(f)\Phi_{w}(f)\overline{\mu(f)}\vert^2}\mu(f)F_z\\
&=\frac{(1-\vert \eta(f)\vert^2)\mu(f)}{\vert1+\eta(f)\overline{\nu(f)} \vert^2-\vert \eta(f)\overline{\mu(f)}\vert^2}F_z\\
&=\frac{(1-\vert \eta(\Psi(F))\vert^2)\mu(\Psi(F))}{\vert1+\eta(\Psi(F))\overline{\nu(\Psi(F))} \vert^2-\vert \eta(\Psi(F))\overline{\mu(\Psi(F))}\vert^2}F_z
\end{align*}
\end{proof}

We now consider the case when the Beltrami equation is of the form $f_{\overline{z}}=\overline{\nu(f) }f_z+\nu(f) \overline{f_z}$. Then \eqref{eq:Quad2} becomes 
\begin{align}\label{eq:Quad3}
\nu(z)+\eta(z)+\overline{\nu(z)}\eta(z)^2=0
\end{align}

Solving the quadratic equation for $\eta$ provided $\nu(z)\neq 0$ gives 
\begin{align}\label{eq:Eta}
\eta_\pm(z)&=-\frac{1}{2\overline{\nu(z)}}\pm \sqrt{\frac{1}{4\overline{\nu(z)}^2}-\frac{\nu(z)}{\overline{\nu(z)}}}
\end{align}
using the principal branch of the square root. 

\begin{Lem}
Let $\nu(z)=re^{i\theta}$. Then 
\begin{align*}
\eta_\pm(z)=\left\{
    \begin{array}{ll}
    \displaystyle -\bigg(\frac{1}{2r}\mp\sqrt{\bigg(\bigg(\frac{1}{2r}\bigg)^2-1\bigg)}\bigg)e^{i\theta}, & \text{if $\theta\in [-\pi/2,\pi/2]$}\\
   \displaystyle -\bigg(\frac{1}{2r}\pm \sqrt{\bigg(\bigg(\frac{1}{2r}\bigg)^2-1\bigg)}\bigg)e^{i\theta}& \text{if $\theta\in (-\pi,-\pi/2)\cup(\pi/2,\pi)$}.
\end{array} \right.
\end{align*}
\end{Lem}

\begin{proof}
If $-\frac{\pi}{2}\leq \theta\leq \frac{\pi}{2}$ then 
\begin{align*}
\frac{1}{4\overline{\nu(z)}^2}-\frac{\nu(z)}{\overline{\nu(z)}}&=\frac{1}{4r^2e^{-2i\theta}}-\frac{re^{i\theta}}{re^{-i\theta}}\\
&=\bigg(\bigg(\frac{1}{2r}\bigg)^2-1\bigg)e^{2i\theta}
\end{align*}
and 
\begin{align*}
\eta_\pm(z)=-\frac{1}{2\overline{\nu(z)}}\pm\sqrt{\frac{1}{4\overline{\nu(z)}^2}-\frac{\nu(z)}{\overline{\nu(z)}}}=-\bigg(\frac{1}{2r}\mp\sqrt{\bigg(\bigg(\frac{1}{2r}\bigg)^2-1\bigg)}\bigg)e^{i\theta}
\end{align*}
using the principal argument. If $\frac{\pi}{2}<\theta\leq \pi$ we write $\theta=\frac{\pi}{2}+\phi$ with $0<\phi\leq \frac{\pi}{2}$. Thus 
\begin{align*}
2\theta=\pi+2\phi=-\pi+2\phi
\end{align*}
using the principal argument. This gives
\begin{align*}
\eta_\pm(z)=-\frac{e^{i\pi/2+i\phi}}{2r}\pm e^{-i\pi/2+i\phi}\sqrt{\bigg(\bigg(\frac{1}{2r}\bigg)^2-1\bigg)}=-\bigg(\frac{1}{2r}\pm \sqrt{\bigg(\bigg(\frac{1}{2r}\bigg)^2-1\bigg)}\bigg)e^{i\pi/2+i\phi}.
\end{align*}
Similarly, if $-\pi<\theta< -\frac{\pi}{2}$ we write $\theta=-\pi/2-\phi$, with $0<\phi\leq \frac{\pi}{2}$. Thus
\begin{align*}
2\theta=-\pi-2\phi=\pi-2\phi
\end{align*}
using the principal argument. 
This gives
\begin{align*}
\eta_\pm(z)=-\frac{e^{-i\pi/2-i\phi}}{2r}\pm e^{i\pi/2-i\phi}\sqrt{\bigg(\bigg(\frac{1}{2r}\bigg)^2-1\bigg)}=-\bigg(\frac{1}{2r}\pm \sqrt{\bigg(\bigg(\frac{1}{2r}\bigg)^2-1\bigg)}\bigg)e^{-i\pi/2-i\phi}.
\end{align*}
\end{proof}

On the other hand whenever $\nu(z)=0$ \eqref{eq:Quad3} implies that $\eta(z)=0$. In order for $\Phi$ to solve a locally uniformly elliptic Beltrami equation on $N$ we want either that $\vert \eta(z)\vert<1$ or $\vert \eta(z)\vert>1$ or locally. Since we may need to choose $\eta(z)=0$ where $\nu(z)=0$, the only root compatible with this condition is the root whose modulus is less than 1. Thus we define
\begin{align}\label{def:Eta}
\eta(z)&=\left\{
    \begin{array}{ll}
    0, & \text{if $\nu(z)=0$}\\
   \eta_-(z),& \text{if $\arg(\nu(z))\in[-\pi/2,\pi/2] $}\\
   \eta_+(z), &\text{if $\arg(\nu(z))\in (-\pi,-\pi/2)\cup(\pi/2,\pi)$}
   \end{array} \right. \nonumber \\
   &=\left\{
    \begin{array}{ll}
    0, & \text{if $\nu(z)=0$}\\
   \displaystyle -\bigg(\frac{1}{2\vert \nu(z)\vert}- \sqrt{\bigg(\frac{1}{2\vert \nu(z)\vert}\bigg)^2-1)}\bigg)e^{i\arg(\nu(z))} &
   \text{otherwise},
   \end{array} \right. \\
&=\left\{
    \begin{array}{ll}
    0, & \text{if $\nu(z)=0$}\\
   \displaystyle -\bigg(\frac{1}{2\vert \nu(z)\vert}- \sqrt{\bigg(\frac{1}{2\vert \nu(z)\vert}\bigg)^2-1}\bigg)\frac{\nu(z)}{\vert \nu(z)\vert} &
   \text{otherwise},
   \end{array} \right. \nonumber\\
&=\left\{
    \begin{array}{ll}
    0, & \text{if $\nu(z)=0$}\\
   \displaystyle -\frac{1}{2\overline{\nu(z)}}+\frac{1}{2\overline{\nu(z)}}\sqrt{1-4\vert \nu(z)\vert^2}&
   \text{otherwise}.
   \end{array} \right. \nonumber
\end{align}

\begin{Lem}[Ellipticity of $\eta$]
Whenever $\vert \nu(z)\vert<1/2$, $\vert \eta(z)\vert<1$.
Moreover as $\eta(z)\to 0$ as $\nu(z)\to 0$. 
\end{Lem}

\begin{proof}
For $r\neq 0$
\begin{align*}
\vert \eta(z)\vert=\frac{1}{2r}-\sqrt{\bigg(\bigg(\frac{1}{2r}\bigg)^2-1\bigg)}.
\end{align*}
In addition,
\begin{align*}
\frac{1}{2r}+\sqrt{\bigg(\bigg(\frac{1}{2r}\bigg)^2-1\bigg)}>\frac{1}{2r}>1.
\end{align*}
Since $\displaystyle\vert \eta_+(z)\vert \vert \eta_-(z)\vert=\bigg\vert \frac{\nu(z)}{\overline{\nu(z)}}\bigg\vert=1$ it follows that $\vert \eta(z)\vert<1$. Since 
\begin{align*}
\lim_{r\to 0^+}\frac{1}{2r}+\sqrt{\bigg(\bigg(\frac{1}{2r}\bigg)^2-1\bigg)}=+\infty
\end{align*}
it follows that $\eta(z)\to 0$ as $\nu(z)\to 0$.
\end{proof}

\begin{Lem}[Existence of $\Phi$]
\label{lem:ex1}
Let $N\subset \C$ be a finitely connected domain conformally equivalent to a bounded Koebe domain. Let $\eta\in L^\infty(N,\C)$ and assume that for every $V\Subset N$ there exists a $k=k(V)$ such that $\Vert \eta\Vert_{L^\infty(V)}=k<1$. Then there exists an orientation preserving homeomorphic solution in $W^{1,2}_{loc}(N,N)$ to
\begin{align*}
\Phi_{\overline{z}}(z)=\eta(z)\Phi_z(z).
\end{align*}
If $\Vert \eta\Vert_{L^\infty(N)}=k<1$ the assumption that $N$ is necessarily conformally equivalent to a \emph{bounded} Koebe domain can be removed, in particular when $N=\C$. 
\end{Lem}

\begin{proof}
The proof is a direct generalisation of the first part of the proof of \cite[Theorem 4.1]{ADPZ}. In the case when $\Phi$ solves a uniformly elliptic equation and $N$ is not necessarily conformally equivalent to a bounded Koebe domain, this follows from the measurable Riemann mapping theorem, see \cite[Theorem 5.3.4]{AIM}. 
\end{proof}

By the theory of quasiregular maps (\cite[Corolllary 3.10.3]{AIM}), it follows that $\Phi$ is locally H÷lder continuous and one verifies that $\Phi$ satisfies the assumptions of Lemma \ref{lem:Chain}, justifying the use of the chain rule.

\begin{rem}
In fact for the application of Lemma \ref{lem:ex1} to Lemma \ref{lem:CB1} any homeomorphism $\Phi$ will do, not necessarily ones such that $\Phi:N\to N$. Moreover, if $N$ is not conformally equivalent to a bounded Koebe domain and the Beltrami equation is degenerate, then \cite[Theorem 4.1]{ADPZ} does not apply and we do not know of a general theorem which guarantees the existence of homeomorphic solutions. However, in some cases one can prove existence directly by explicit methods. This is the case for the minimal surface equation in Example \ref{ex:Mini}. There
\begin{align*}
\Phi(z)=\frac{2z}{1+\sqrt{1+4\vert z\vert^2}}
\end{align*}
is a homeomorphism $\Phi: \C\to \Di$ which solves the degenerate Beltrami equation
\begin{align*}
\Phi_{\overline{z}}(z)=-\frac{2\vert z\vert^2}{2\vert z\vert^2+1+\sqrt{1+4\vert z\vert^2}}\frac{z}{\overline{z}}\Phi_z(z),
\end{align*}
since $\vert \frac{2\vert z\vert^2}{2\vert z\vert^2+1+\sqrt{1+4\vert z\vert^2}}\vert \to 1$ as $\vert z\vert\to +\infty.$
\end{rem}

\begin{rem}\label{rem:1}
If we had chosen the rot such that $\vert \eta(z)\vert >1$ whenever $\nu(z)\neq 0$ we would still have to choose $\eta(z)=0$ whenever $\nu(z)=0$. If $\nu$ is such that there exists an open set $U\subset N$, with $U\neq N$ and $\nu(z)=0$ for $z\in U$ but $\nu(z)\neq 0$ for a.e $z\in N\setminus U$ then, there can exists no homeomorphism $\Phi$ solving \eqref{eq:Ansatz}. Indeed, by \cite[Theorem 3.3.5]{AIM} for any homeomorphism $\Phi\in W^{1,1}_{loc}(\Om,\C)$, the Jacobian $J(z,\Phi)=\det(D\Phi(z))$ does not change sign if $\Phi$ is a homeomorphism. On the other hand if $\Phi$ solves \eqref{eq:Ansatz} with $\eta$ chosen as described, then 
\begin{align*}
J(z,\Phi)=\vert \Phi_z(z)\vert^2-\vert \Phi_{\overline{z}}(z)\vert^2
=\left\{
    \begin{array}{ll}
    \vert \Phi_z(z)\vert^2>0 & \text{ for a.e. $z\in U$},\\
     (1-\vert \nu_-(z)\vert^2)\vert \Phi_z(z)\vert^2<0& \text{ for a.e. $z\in N\setminus U$},
\end{array} \right.
\end{align*}
a contradiction, since $\Phi_z(z)\neq 0$ for a.e. $z\in N$. Even though the approach in \cite{KP} is different from ours, they choose an orientation reversing map in \cite{KP}. Since this map is unique up to composition with holomorphic maps, this means that they have chosen $\nu$ such that $\vert \nu(z)\vert>1$ whenever $\nu(z)=0$ instead. However, if $\A$ is merely smooth, rather than real analytic we may have $\nu(z)=0$ on a set of positive measure, which means that we cannot chose $\eta$ such that $\vert \eta(z)\vert >1$ for $a.e.$ $z$ in this case. This possibility does not seem to have been considered in \cite{KP}. Indeed, in section \ref{sec:CounterEx} we give an example of structure field $\A$ where this happens. 
\end{rem}

The equation for $F$ in this case becomes
\begin{align*}
F_{\overline{z}}=\frac{(1-\vert \eta(\Psi(F))\vert^2)\overline{\nu(\Psi(F))}}{\vert1+\eta(\Psi(F))\overline{\nu(\Psi(F))} \vert^2-\vert \eta(\Psi(F))\nu(\Psi(F))\vert^2}F_z
\end{align*}
which holds for a.e. $z\in \Om$ provided the denominator is nonzero. 

\begin{Def}
We define for $z\in N$
\begin{align}\label{def:gamma}
\gamma(z)=\frac{ (1-\vert\eta(z)\vert^2)\overline{\nu(z)}}{\vert1+\eta(z)\overline{\nu(z)} \vert^2-\vert \eta(z)\nu(z)\vert^2} 
\end{align}
\end{Def}

\begin{Lem}[Ellipticity of $\gamma$]
\label{lem:Egamma}
$\vert1+\eta(z)\overline{\nu(z)} \vert^2-\vert \eta(z)\nu(z)\vert^2\neq 0$ for a.e $z\in N$. Moreover, $\vert \gamma(z)\vert<1$ whenever $2\vert \nu(z)\vert<1$. Furthermore, $\gamma(z)=-\overline{\eta(z)}$.
\end{Lem}

\begin{proof}
Using that same $\nu(z)=re^{i\theta}$ and $\displaystyle \eta(z)=-\bigg(\frac{1}{2r}-\sqrt{\bigg(\bigg(\frac{1}{2r}\bigg)^2-1\bigg)}\bigg)e^{i\theta}$ we get 
that $\eta(z)\overline{\nu(z)}\in \R$ and $\eta(z)\overline{\nu(z)}=-\vert \eta(z)\nu(z)\vert$. 
\begin{align*}
&\vert 1+\eta(z)\overline{\nu(z)} \vert^2-\vert \eta(z)\nu(z)\vert^2=1+2\eta(z)\overline{\nu(z)}\\
&=1-2r\bigg(\frac{1}{2r}-\sqrt{\bigg(\bigg(\frac{1}{2r}\bigg)^2-1\bigg)}\bigg)\\
&=2r\sqrt{\bigg(\bigg(\frac{1}{2r}\bigg)^2-1\bigg)}=\sqrt{1-(2r)^2}>0
\end{align*}
for all $r<1/2$. Thus $\vert1+\eta(z)\overline{\nu(z)} \vert^2-\vert \eta(z)\nu(z)\vert^2\neq 0$ for a.e $z\in N$.
Furthermore,
\begin{align*}
&(1-\vert\eta(z)\vert^2)\vert \nu(z)\vert=r\bigg(1-\bigg(\frac{1}{2r}-\sqrt{\bigg(\bigg(\frac{1}{2r}\bigg)^2-1\bigg)}\bigg)^2\bigg)\\
&=r\bigg(1-\bigg(\frac{1}{2r}\bigg)^2+\frac{1}{r}\sqrt{\bigg(\bigg(\frac{1}{2r}\bigg)^2-1\bigg)}-\bigg(\frac{1}{2r}\bigg)^2+1\bigg)\\
&=r\bigg(2-2\bigg(\frac{1}{2r}\bigg)^2+\frac{1}{r}\sqrt{\bigg(\bigg(\frac{1}{2r}\bigg)^2-1\bigg)}\bigg).
\end{align*}
If we let $\displaystyle p=\sqrt{\bigg(\bigg(\frac{1}{2r}\bigg)^2-1\bigg)}$ then we can write the above as 
\begin{align*}
&(1-\vert\eta(z)\vert^2)\vert \nu(z)\vert=-2rp^2+p.
\end{align*}
Thus 
\begin{align}
\vert \gamma(z)\vert=\frac{-2rp^2+p}{-2rp}=-\frac{1}{2r}+\sqrt{\bigg(\bigg(\frac{1}{2r}\bigg)^2-1\bigg)}
\end{align}
Thus $\vert \gamma(z)\vert=\vert \eta(z)\vert$ and so $\vert \gamma(z)\vert<1$ whenever $\vert \nu(z)\vert<1/2$. In addition 
\begin{align*}
\gamma(z)=-(\frac{1}{2r}-\sqrt{\bigg(\bigg(\frac{1}{2r}\bigg)^2-1\bigg)}\bigg)e^{-i\theta}=-\overline{\eta(z)}. 
\end{align*}
\end{proof}

\begin{rem}
By Lemma \ref{lem:Egamma} $\Vert \mu(w)\Vert_{L^{\infty}(U)}=k_U<1$ for all $U\Subset N$. 
\end{rem}

Thus, to conclude the complex valued field $F=\Phi\circ f=\Phi\circ u_z$ solves the $\C$-quasilinear Beltrami equation 
\begin{align}\label{eq:CB}
F_{\overline{z}}(z)=\mu(\Phi^{-1}(F(z)))F_z(z)
\end{align}
for a.e. $z\in \Om$ and where $\mu$ is given as in Definition \ref{def:ConfB}. 

This concludes the proof of Theorem \ref{thm:Main}.


\subsection{\sffamily The anti-$\C$-linear Beltrami equation and pseudo-analytic functions}\label{sec:AntiC}

Consider the anti-$\C$-linear Beltrami equation

\begin{align}\label{eq:AntiC}
g_{\overline{z}}(z)=-\phi(z)\overline{g_z(z)}, \quad z\in U 
\end{align}
where $\Vert \phi\Vert_{L^\infty(V)}=k_V<1$ for all $V\Subset U$. Dimer models (see \cite{ADPZ}) have the special property that the Beltrami coefficient $\phi$ is a holomorphic function. This makes equation \eqref{eq:AntiC} exact in the sense 
\begin{align*}
g_{\overline{z}}(z)+\phi(z)\overline{g_z(z)}=(g(z)+\phi(z)\overline{g(z)})_{\overline{z}}=0,
\end{align*}
and so the function $h=g+\phi\overline{g}$ is holomorphic. This allows one to parametrise solutions of \eqref{eq:AntiC} by means of holomorphic functions. We could ask if something similar is true for more general anti-$\C$-linear Beltrami equation where $\phi$ is no longer holomorphic. If $\phi_{\overline{z}}$ exists in a distributional sense as an $L^1$ function and $\phi_{\overline{z}}g,\phi g_z\in L^1_{loc}(U)$, the Leibniz rule implies
\begin{align*}
(-\phi(z)\overline{g(z)})_{\overline{z}}=-\phi_{\overline{z}}(z)\overline{g(z)}-\phi(z)\overline{g_z(z)}
\end{align*}
and equation \eqref{eq:AntiC} can be written as 
\begin{align*}
g_{\overline{z}}(z)=(-\phi(z)\overline{g(z)})_{\overline{z}}+\phi_{\overline{z}}(z)\overline{g(z)}
\end{align*}

or equivalently 

\begin{align*}
(g(z)+\phi(z)\overline{g(z)})_{\overline{z}}=\phi_{\overline{z}}(z)\overline{g(z)}
\end{align*}

Define the function 
\begin{align*}
h(z):=g(z)+\phi(z)\overline{g(z)}. 
\end{align*}

Then elementary algebra yields
\begin{align*}
g(z)=\frac{1}{1-\vert \phi(z)\vert^2}(h(z)-\phi(z)\overline{h(z)})
\end{align*}

and $h$ solves the inhomogeneous linear $\overline{\dv}$-equation
\begin{align*}
h_{\overline{z}}(z)=\frac{\phi_{\overline{z}}(z)}{1-\vert \phi(z)\vert^2}(\overline{h(z)}-\overline{\phi(z)}h(z))
\end{align*}

Let $\Om\subset \C$ be a bounded domain. A function $h$ which satisfy the equation 
\begin{align}\label{eq:GHolo}
h_{\overline{z}}(z)+\alpha(z)h(z)+\beta(z)\overline{h(z)}=0
\end{align}
in $\Om$ for some $\alpha,\beta\in L^{p}(\Om,\C)$, $p>2$ is a \emph{generalized analytic function in the sense of Vekua}, see the monograph \cite{V}. In particular, \cite[ch. 3, p. 144-146]{V} shows that every solution of \eqref{eq:GHolo} belongs to $W^{1,p}(\Om,\C)$. Moreover every solution $h$ is of the form 
\begin{align*}
h(z)=\psi(z)e^{\omega(z)},
\end{align*}
where 
\begin{align}\label{eq:Omega}
\omega(z)=\frac{1}{\pi}\int_{\Om}\frac{\alpha(z)+\beta(z)(\chi_{h(z)=0}+\frac{\overline{h(z)}}{h(z)}\chi_{h(z)\neq 0})}{w-z}dA(w)
\end{align}
and $\psi$ is a holomorphic function on $\Om$. Thus we have shown that if $\Vert \phi\Vert_{L^\infty(V)}=k_V<1$ for all $V\Subset U$ and in addition $\phi\in W^{1,p}(U)$, then every solution $g$ of \eqref{eq:AntiC} is of the form 
\begin{align*}
g(z)=\frac{1}{1-\vert \phi(z)\vert^2}(h(z)-\phi(z)\overline{h(z)}).
\end{align*}
where $h$ is a solution of \eqref{eq:GHolo}. This completes the proof of Theorem \ref{Thm:Main2}.


\subsection{\sffamily Comparison to other types of linearisations}

We will now conclude this section by discussing other types of linearisations and how they differ from the aforementioned one. These methods will not be employed in this paper but merely serve as a comparison. Rather than considering the complex gradient and applying the chain rule to derive the equation \eqref{eq:Cgrad} from the \eqref{eq:LL}, we could instead consider the $\A$-harmonic conjugate $v$ of $u$ defined according to 
\begin{align*}
\nabla v(z)=\star \A(\nabla u(z)),
\end{align*}
where the Hodge star $\star$ is given by multiplication by the matrix
\begin{align*}
\begin{bmatrix}
0 & -1\\
1 & 0
\end{bmatrix}.
\end{align*}

Such $v$ always exists and is unique up to constant on a simply connected domain. Defining the complex valued field $\mathscr{F}=u+iv$ and the considering the nonlinear Cayley transform 
\begin{align*}
\mathcal{H}(w)=(I-\A)\circ(I+\A)^{-1}(\overline{w})
\end{align*}
as in \cite{ACFJK17,ACFJK20,ADPZ}, one can show that $\mathscr{F}$ solves the fully nonlinear Beltrami equation 
\begin{align}\label{eq:FLB}
\mathscr{F}_{\overline{z}}(z)=\mathcal{H}(\mathscr{F}_z(z)).
\end{align}
From here one can do an Iwaniec-Sbordonne linearisation as in \cite[ch.16]{AIM} as follows: Define 
\begin{align*}
k(z)=\frac{\vert \mathscr{F}_{\overline{z}}(z)\vert}{\vert \mathscr{F}_z(z)\vert}
\end{align*}
and let 
\begin{align*}
\mathbf{n}(z)=\frac{\mathscr{F}_{\overline{z}}(z)-k(z)\mathscr{F}_z(z)}{\vert \mathscr{F}_{\overline{z}}(z)-k(z)\mathscr{F}_z(z)\vert}
\end{align*}
if $\mathscr{F}_{\overline{z}}(z)-k(z)\mathscr{F}_z(z)\neq 0$, and otherwise let $n(z)$ be a unit vector orthogonal to both $\mathscr{F}_{\overline{z}}(z)$ and $\mathscr{F}_z(z)$. Define the measurable linear transformation $\mathcal{M}(z)$ to be  
\begin{align*}
\mathcal{M}(z)=k(z)[I-2\mathbf{n}(z)\otimes \mathbf{n}(z)].
\end{align*}
The any solution to \eqref{eq:FLB} also solves the linear equation 
\begin{align}\label{eq:FLB2}
\mathscr{F}_{\overline{z}}(z)=\mathcal{M}(z)\mathscr{F}_z(z). 
\end{align}
Unfortunately the regularity of the coefficients of the linearised equation \emph{depends on $\mathscr{F}$} itself. Furthermore, it also fails if $\mathcal{H}$ is not a $k$-Lipschitz function for some $k<1$ as the resulting linear equation becomes degenerate. If one does not have an a priori estimate of $\mathscr{F}_z$ itself, which is what one wants to achieve one has no way of controlling the degeneracy of the resulting linear equation. This is the case for example when \eqref{eq:LL} is the $p$-Laplace equation. In the case when $p=3$, $\mathcal{H}$ can be explicitly computed to give
\begin{align*}
\mathcal{H}(w)=\frac{2(1+\sqrt{1+4\vert w\vert}-2\vert w\vert)}{(1+\sqrt{1+4 \vert w\vert})^2}\overline{w}
\end{align*}

Differentiating \eqref{eq:FLB} with respect to $\dv_z$ and setting $\mathfrak{f}=\mathscr{F}_z$ give the quasilinear equation 
\begin{align*}
\mathfrak{f}_{\overline{z}}(z)=\frac{\mathcal{H}_w(\mathfrak{f})}{1-\vert \mathcal{H}_{\overline{w}}(\mathfrak{f})\vert^2(\mathfrak{f})}\mathfrak{f}_z(z)+\frac{\overline{\mathcal{H}_w(\mathfrak{f})}\mathcal{H}_{\overline{w}}(\mathfrak{f})}{1-\vert \mathcal{H}_{\overline{w}}(\mathfrak{f})\vert^2}\overline{\mathfrak{f}_z(z)},
\end{align*}
which is \emph{not uniformly elliptic} in the case of the $p$-Laplacian despite the fact that the complex gradient solves the uniformly elliptic Beltrami equation \eqref{eq:Cgrad}.


\section{\sffamily Applications to highly degenerate elliptic equations}

\begin{ex}
This example is from \cite{BS} and arises in the study of maximal spacelike hypersurfaces in Minkowski space. Let $\Om\subset \R^n$ be a domain and let $\mathscr{A}=\{u\in W^{1,\infty}(\Om): \nabla u(x)\in B_1(0)\}$ denote the space of admissible Lipschitz functions on $\Om$ whose gradient lie in the unit ball. Consider the functional 
\begin{align*}
I[u]=\int_\Om-\sqrt{1-\vert \nabla u(x)\vert^2}dx
\end{align*}
defined on $\mathscr{A}$. We want to study the associated Euler-Lagrange equations in the case when $\Om \subset \C$ is a domain in the plane and $B_1(0)=\Di$. We have  
\begin{align*}
\nabla \sigma(z)=\A(z)=\frac{z}{\sqrt{1-\vert z\vert^2}}=\frac{\vert z\vert}{\sqrt{1-\vert z\vert^2}}\frac{z}{\vert z\vert}.
\end{align*}
This gives

\begin{align*}
\mathbf{A}(z)=\frac{2\vert z\vert}{\sqrt{1-4\vert z\vert^2}}\frac{z}{\vert z\vert}
\end{align*}
and $N=1/2\Di$. Let 
\begin{align*}
\rho(t)=\frac{2t}{\sqrt{1-4t^2}}, \quad \dot{\rho}(t)=\frac{2}{(1-4t^2)^{3/2}}.
\end{align*}
so that $\mathbf{A}(z)=\rho(\vert z\vert)z/\vert z\vert$. Then 
\begin{align*}
\mathbf{A}_z(z)&=\frac{1}{2}\frac{4-8\vert z\vert^2}{(1-4\vert z\vert^2)^{3/2}}\\
\mathbf{A}_{\overline{z}}(z)&=-\frac{1}{2}\frac{z}{\overline{z}}\frac{8\vert z\vert^2}{(1-4\vert z\vert^2)^{3/2}}.
\end{align*}
Thus

\begin{align*}
\nu(z)=\frac{1}{2}\frac{z}{\overline{z}}\frac{8\vert z\vert^2}{4-8\vert z\vert^2}=\frac{z}{\overline{z}}\frac{2\vert z\vert^2}{1-2\vert z\vert^2}
\end{align*}

This gives 
\begin{align*}
\eta(z)=-\frac{1-2\vert z\vert^2-\sqrt{1-4\vert z\vert^2}}{2\vert z\vert^2}\frac{z}{\overline{z}}
\end{align*}
and hence the conformal structure of the Euler-Lagrange equation is 
\begin{align*}
-\frac{1-2\vert u_z\vert^2-\sqrt{1-4\vert u_z\vert^2}}{2\vert u_z\vert^2}\frac{u_z}{\overline{u_z}}.
\end{align*}

We want to find a homeomorphic solution $\Phi: \frac{1}{2}\Di\to \frac{1}{2}\Di$ of $\Phi_{\overline{z}}=\eta \Phi_z$. Exploiting the fact that $\eta$ is of the form $\eta(z)=\gamma(\vert z\vert)\frac{z}{\overline{z}}$ where 
\begin{align*}
\gamma(t)=-\frac{1-2t^2-\sqrt{1-4t^2}}{2t^2}
\end{align*}
we make the ansatz that $\Phi$ is a radial stretching map, i.e., $\Phi(z)=\rho(\vert z\vert)\frac{z}{\vert z\vert}$. This gives 

\begin{align*}
\Phi_{\overline{z}}(z)&=\frac{1}{2}\frac{z}{\overline{z}}\bigg[\dot{\rho}(\vert z\vert)-\frac{\rho(\vert z\vert)}{\vert z\vert}\bigg]\\
\Phi_z(z)&=\frac{1}{2}\bigg[\dot{\rho}(\vert z\vert)+\frac{\rho(\vert z\vert)}{\vert z\vert}\bigg]
\end{align*}
which implies that
\begin{align*}
\frac{1}{2}\frac{z}{\overline{z}}\bigg[\dot{\rho}(\vert z\vert)-\frac{\rho(\vert z\vert)}{\vert z\vert}\bigg]=\frac{1}{2}\gamma(\vert z\vert)\frac{z}{\overline{z}}\bigg[\dot{\rho}(\vert z\vert)+\frac{\rho(\vert z\vert)}{\vert z\vert}\bigg]
\end{align*}

which reduces to the ODE
\begin{align*}
\dot{\rho}(t)=\frac{1}{t}\frac{1+\gamma(t)}{1-\gamma(t)}\rho(t). 
\end{align*}
This ODE is separable and we find that

\begin{align*}
\ln \vert \rho(t)\vert=\int \frac{1}{t}\frac{\gamma(t)+1}{\gamma(t)-1}dt=\int \frac{1}{t}+\frac{4t}{\sqrt{1-4t^2}-1} dt=\log\vert t\vert-\sqrt{1-4t^2}-\log\vert 1-\sqrt{1-4t^2}\vert+C
\end{align*}

Thus, the solution $\Phi(z)$ is given by 
\begin{align*}
\Phi(z)=\frac{e^C\vert z\vert e^{-\sqrt{1-4\vert z\vert^2}}}{1-\sqrt{1-4\vert z\vert^2}}\frac{z}{\vert z\vert}=\frac{e^Ce^{-\sqrt{1-4\vert z\vert^2}}}{1-\sqrt{1-4\vert z\vert^2}}z
\end{align*}

Choosing $C=\log(1/2)$ gives a homeomorphism $\Phi:\frac{1}{2}\Di\to \frac{1}{2}\Di$.  We cannot invert $\Phi$ explicitly in terms of elementary functions, however, all regularity of $\Phi^{-1}$ can be deduced. 
To conclude, $F=\Phi(u_z)$ solves the $\C$-quasilinear Beltrami equation 
\begin{align*}
F_{\overline{z}}(z)=-\overline{\eta(\Phi^{-1}(F(z))}F_z(z).
\end{align*}
\end{ex}

\begin{ex}
We again consider the example of the $p$-orthotropic Laplacian in the plane. 
Let $\Om \subset \C$ be a bounded domain and consider the $p$-orthotropic functional
\begin{align}\label{eq:PO}
I[u]=\int_{\Om}(\vert u_x\vert^p+\vert u_y\vert^p)dxdy
\end{align}
defined on $W^{1,p}(\Om)$. Any critical point of \eqref{eq:PO} is a weak solution of the Euler-Lagrange equations
\begin{align}\label{eq:OPL}
(\vert u_x\vert^{p-2}u_x)_x+(\vert u_y\vert^{p-2}u_y)_y=0. 
\end{align}
By \cite{BS} every weak solution is in $C^1(\Om)$. With $\sigma(x,y)=\vert x\vert^p+\vert y\vert^p$ we have 
\begin{align*}
\A(z)=\nabla \sigma(z)=2^{1-p}(\vert z+\overline{z}\vert^{p-2}(z+\overline{z})+\vert z-\overline{z}\vert^{p-2}(z-\overline{z}))
\end{align*}
where $z=x+iy$. Thus
\begin{align*}
\mathbf{A}(z)=\overline{\A(2\overline{z})}=\vert z+\overline{z}\vert^{p-2}(z+\overline{z})+\vert z-\overline{z}\vert^{p-2}(z-\overline{z})
\end{align*}

We see that $\mathbf{A}(z)$ is the sum of the compositions of the radial stretching map
\begin{align*}
R(z)=\rho(\vert z\vert)\frac{z}{\vert z\vert}
\end{align*}
with $\rho(\vert z\vert)=\vert z\vert^{p-1}$ with the maps $z\mapsto z+\overline{z}$ and $z\mapsto z-\overline{z}$. Using the chain rule and the formulas 

\begin{align*}
R_z(z)&=\frac{1}{2}\bigg[\dot{\rho}(\vert z\vert)+\frac{\rho(\vert z\vert)}{\vert z\vert}\bigg]\\
R_{\overline{z}}(z)&=\frac{1}{2}\frac{z}{\overline{z}}\bigg[\dot{\rho}(\vert z\vert)-\frac{\rho(\vert z\vert)}{\vert z\vert}\bigg]\
\end{align*}

for the radial stretching maps we get 
\begin{align*}
\mathbf{A}_z(z)&=R_z(z+\overline{z})+R_{\overline{z}}(z+\overline{z})+R_z(z-\overline{z})-R_{\overline{z}}(z-\overline{z})\\
&=\frac{1}{2}\bigg[\dot{\rho}(z+\overline{z})+\frac{\rho(z+\overline{z})}{\vert z+\overline{z}\vert}\bigg]+\frac{1}{2}\bigg[\dot{\rho}(z+\overline{z})-\frac{\rho(z+\overline{z})}{\vert z+\overline{z}\vert}\bigg]\frac{z+\overline{z}}{\overline{z+\overline{z}}}\\
&+\frac{1}{2}\bigg[\dot{\rho}(z-\overline{z})+\frac{\rho(z-\overline{z})}{\vert z-\overline{z}\vert}\bigg]-\frac{1}{2}\bigg[\dot{\rho}(z-\overline{z})-\frac{\rho(z-\overline{z})}{\vert z-\overline{z}\vert}\bigg]\frac{z-\overline{z}}{\overline{z-\overline{z}}}\\
&=\frac{1}{2}\bigg[\dot{\rho}(z+\overline{z})+\frac{\rho(z+\overline{z})}{\vert z+\overline{z}\vert}\bigg]+\frac{1}{2}\bigg[\dot{\rho}(z+\overline{z})-\frac{\rho(z+\overline{z})}{\vert z+\overline{z}\vert}\bigg]\\
&+\frac{1}{2}\bigg[\dot{\rho}(z-\overline{z})+\frac{\rho(z-\overline{z})}{\vert z-\overline{z}\vert}\bigg]+\frac{1}{2}\bigg[\dot{\rho}(z-\overline{z})+\frac{\rho(z-\overline{z})}{\vert z-\overline{z}\vert}\bigg]\\
&=\dot{\rho}(z+\overline{z})+\dot{\rho}(z-\overline{z})
\end{align*}
and 
\begin{align*}
\mathbf{A}_{\overline{z}}(z)&=R_z(z+\overline{z})+R_{\overline{z}}(z+\overline{z})-R_z(z-\overline{z})+R_{\overline{z}}(z-\overline{z})\\
&=\frac{1}{2}\bigg[\dot{\rho}(z+\overline{z})+\frac{\rho(z+\overline{z})}{\vert z+\overline{z}\vert}\bigg]+\frac{1}{2}\bigg[\dot{\rho}(z+\overline{z})-\frac{\rho(z+\overline{z})}{\vert z+\overline{z}\vert}\bigg]\frac{z+\overline{z}}{\overline{z+\overline{z}}}\\
&-\frac{1}{2}\bigg[\dot{\rho}(z-\overline{z})+\frac{\rho(z-\overline{z})}{\vert z-\overline{z}\vert}\bigg]+\frac{1}{2}\bigg[\dot{\rho}(z-\overline{z})-\frac{\rho(z-\overline{z})}{\vert z-\overline{z}\vert}\bigg]\frac{z-\overline{z}}{\overline{z-\overline{z}}}\\
&=\frac{1}{2}\bigg[\dot{\rho}(z+\overline{z})+\frac{\rho(z+\overline{z})}{\vert z+\overline{z}\vert}\bigg]+\frac{1}{2}\bigg[\dot{\rho}(z+\overline{z})-\frac{\rho(z+\overline{z})}{\vert z+\overline{z}\vert}\bigg]\\
&-\frac{1}{2}\bigg[\dot{\rho}(z-\overline{z})+\frac{\rho(z-\overline{z})}{\vert z-\overline{z}\vert}\bigg]-\frac{1}{2}\bigg[\dot{\rho}(z-\overline{z})-\frac{\rho(z-\overline{z})}{\vert z-\overline{z}\vert}\bigg]\\
&=\dot{\rho}(z+\overline{z})-\dot{\rho}(z-\overline{z})
\end{align*}

We observe that both $\mathbf{A}_z$ and $\mathbf{A}_{\overline{z}}$ are real valued. Hence $\nu(z)$ becomes 
\begin{align*}
\nu(z)=-\frac{\mathbf{A}_{\overline{z}}(z)}{2\mathfrak{R}e[\mathbf{A}_z(z)]}=-\frac{\dot{\rho}(\vert z+\overline{z\vert })-\dot{\rho}(\vert z-\overline{z}\vert)}{2(\dot{\rho}(\vert z+\overline{z}\vert)+\dot{\rho}(\vert z-\overline{z}\vert))}
\end{align*}

Moreover $\dot{\rho}(\vert z\vert)=(p-1)\vert z\vert^{p-2}$. This gives 
\begin{align*}
\nu(z)=-\frac{1}{2}\frac{\vert z+\overline{z}\vert^{p-2} -\vert z-\overline{z}\vert^{p-2}}{\vert z+\overline{z}\vert^{p-2}+\vert z-\overline{z}\vert^{p-2}}
\end{align*}
which is real valued. We see that $\vert \nu(z)\vert=\frac{1}{2}$ whenever either $\mathfrak{R}e[z]=0$ or $\mathfrak{I}m[z]=0$. Thus $\mathbf{A}:\C\to \C$ is \emph{not} $\delta$-monotone. However $\mathbf{A}$ is locally $\delta$-monotone in any connected component of the set $\{z\in \C: \mathfrak{R}e[z]\neq 0 \text{ and }\mathfrak{I}m[z]\neq 0\}$. Let 
\begin{align*}
N=Q_1:=\{z\in \C: \mathfrak{R}e[z]>0 \text{ and }\mathfrak{I}m[z]>0\}
\end{align*}
be the first quadrant which is conformally equivalent to the unit disc. In view of the $C^1$-regularity result in \cite{BS}, the set $U_1=\{z\in \Om: u_z\in Q_1\}$ is open. Thus the equation \eqref{eq:CB} is uniformly elliptic on any open set $V\Subset U$. Moreover, the equation for $\eta$ becomes 
\begin{align*}
\eta(z)^2+\frac{1}{\nu(z)}\eta(z)+1=0
\end{align*}
and we find that 
\begin{align*}
\eta(z)&=\frac{\vert x\vert^{p-2} -\vert y\vert^{p-2}}{\vert x\vert^{p-2}+\vert y\vert^{p-2}-2\sqrt{\vert x\vert^{p-2}\vert y\vert^{p-2}}}
\end{align*}
where $z=x+iy$. In particular $\eta(z)=0$ when $x=y$ and $z\in Q_1$. In addition $\eta(z)=-1$ for $x=0$ and $\eta(z)=1$ for $y=0$ and $x,y>0$.

We would like to study the interior regularity of a weak solution $u$ of \eqref{eq:OPL} up to the boundary of $\Gamma_1=\dv U_1\cap \Om$ as well as the regularity of $\Gamma_1$ itself. Assume that $z_0\in \Gamma_1$ and that $B_r(z_0)\cap U_1$ is simply connected for some $r>0$. We can now first apply Theorem \ref{thm:Main} to solutions of \eqref{eq:OPL} on the open set simply connected set $V=B_r(z_0)\cap U_1$. This shows that $F=\Phi(u_z)$ solves the $\C$-quasilinear equation 
\begin{align}\label{eq:OP2}
F_{\overline{z}}(z)=-\overline{\eta(\Phi^{-1}(F(z))}F_z(z),
\end{align}
where $\Phi:Q_1\to Q_1$ is a homeomorphic solution to 
\begin{align*}
\Phi_{\overline{z}}(z)=\eta(z)\Phi_z(z). 
\end{align*}

We now linearise \eqref{eq:OP2} using Theorem \ref{thm:Main}. Let $\psi: V\to \Di$ be a homeomorphic solution of 
\begin{align*}
\psi_{\overline{z}}(z)=-\overline{\eta(\Phi^{-1}(F(z))}\psi_z(z)
\end{align*}
and set $g=\psi^{-1}$. Then $g$ solves the linear equation 
\begin{align*}
g_{\overline{z}}(f)=\overline{\eta(\Phi^{-1}(h(z)))}\overline{g_z(z)}. 
\end{align*}
where $h: \Di \to Q_1$ is holomorphic. To use these methods to actually deduce regularity of the $p$-orthotropic Laplacian is something that will be done in future work. 
\end{ex}


\section{\sffamily Counterexample}\label{sec:CounterEx}

In section \ref{sec:Proof}, when we consider the reduction to a $\C$-quasilinear equation we consider solutions $\eta$ of the algebraic equation 
\begin{align*}
\nu(z)+\eta(z)+\overline{\nu(z)}\eta(z)^2=0.
\end{align*}
We have two choices of roots, $\eta_1$ and $\eta_2$ such that $\vert \eta_1(z)\vert<1$ and $\vert \eta_2(z)\vert>1$ locally, unless $\nu(z)\neq 0$, in which case $\eta_1(z)=\eta_2(z)=0$. If $\eta$ is real analytic then, then $\nu(z)=0$ on a null set and both equations 
\begin{align}\label{eq:2}
\Phi_{\overline{z}}(z)=\eta_1(z)\Phi_z(z), \quad \Psi_{\overline{z}}(z)=\eta_2(z)\Psi_z(z)
\end{align}
are locally uniformly elliptic, and solutions to the first equation are \emph{orientation preserving} and solutions to the second equation are \emph{orientation reversing} and we could choose either one of them for the reduction to a $\C$-quasilinear Beltrami equation. If on the other hand $\nu$ is merely smooth or worse, then $\nu(z)=0$ need not be a null set. Then the second equation is no longer elliptic. 

Moreover, there exists no such homeomorphic solution $\Psi\in W^{1,1}_{loc}(\Di)$. Indeed, by \cite[Theorem 3.3.5 ]{AIM} the Jacobian $J(z,\Psi)=\det(D\Psi(z))$ does not change sign if $\Psi$ is a homeomorphism. On the other hand if $\Psi$ solves Beltrami equation then
\begin{align*}
J(z,\Psi)=\vert \Psi_z(z)\vert^2-\vert \Psi_{\overline{z}}(z)\vert^2
=\left\{
    \begin{array}{ll}
    \vert \Psi_z(z)\vert^2>0 & \text{ for a.e. }\vert z\vert< 1,\\
     (1-\vert \eta_2(z)\vert^2)\vert \Psi_z(z)\vert^2<0& \text{ for a.e. }\vert z\vert>1,
\end{array} \right.
\end{align*}
since $\Psi_z(z)\neq 0$ for a.e. $z\in \Di$.  

We will now construct a $\delta$-monotone field for which $\nu(z)=0$ on a set of positive measure. Let
\begin{align*}
\A(z)= \left\{
    \begin{array}{ll}
    z, & \vert z\vert\leq 1,\\
    \displaystyle z+\eps\frac{1-\vert z\vert}{1+\vert z\vert}\frac{z}{\vert z\vert}, & \vert z\vert>1.
\end{array} \right.
\end{align*}
where $\eps>0$ is to be determined.
One verifies that $\A:\C\to \C$ is Lipschitz continuous. Indeed, let $w\in \Di$ and $z\in \Di^c$. Then

\begin{align*}
\vert \A(z)-\A(w) \vert=\bigg\vert z+\eps\frac{1-\vert z\vert}{1+\vert z\vert}\frac{z}{\vert z\vert}-w\bigg\vert\leq \vert z-w\vert+\eps\frac{\vert z\vert-1}{\vert z\vert+1},
\end{align*}

and so 
\begin{align*}
\frac{\vert \A(z)-\A(w) \vert}{\vert z-w\vert}&\leq1+\frac{\eps}{\vert z\vert+1}\frac{\vert z\vert-1}{\vert z-w\vert }\leq1+\frac{\eps}{\vert z\vert+1}\frac{\vert z\vert-\vert w\vert}{\vert z-w\vert }\\
&\leq1+\frac{\eps}{\vert z\vert+1}\frac{\vert z\vert-\vert w\vert}{\vert z\vert-\vert w\vert }\leq 1+\frac{\eps}{2}. 
\end{align*}

Moreover, a computation yields

\begin{align*}
\A_z(z)= \left\{
    \begin{array}{ll}
    1, & \vert z\vert< 1,\\
    \displaystyle 1+\frac{\eps}{2}\frac{1-2\vert z\vert-\vert z\vert^2}{\vert z\vert(1+\vert z\vert)^2}, & \vert z\vert>1.
\end{array} \right.
\end{align*}

\begin{align*}
\A_{\overline{z}}(z)= \left\{
    \begin{array}{ll}
    0, & \vert z\vert< 1,\\
    \displaystyle \frac{\eps}{2}\frac{z}{\overline{z}}\frac{\vert z\vert^2-2\vert z\vert-1}{\vert z\vert(1+\vert z\vert)^2}, & \vert z\vert>1.
\end{array} \right.
\end{align*}

One notes that $\A_z(z)$ is real. By choosing $\eps>0$ sufficiently small we can ensure that $\A_z(z)>\frac{1}{2}$ for a.e $z\in 2\Di$. Moreover, for $\eps>0$ sufficiently small we can make $\vert \A_{\overline{z}}(z)\vert <1/4$. Thus $\A$ satisfies \eqref{eq:dm} with $\delta=\sqrt{3}/2$. In addition, $\mathbf{A}_z(z)=2\overline{\A_w(2\overline{z})}$ and $\mathbf{A}_z(z)=2\overline{\A_{\overline{w}}(2\overline{z})}$. This gives  
\begin{align*}
\nu(z)= \left\{
    \begin{array}{ll}
    0, & \vert z\vert< 1/2,\\
    \displaystyle \frac{1}{2}\frac{z}{\overline{z}}\frac{\vert 2z\vert^2-2\vert 2z\vert-1}{2\vert 2z\vert(1+\vert 2z\vert)^2+\eps(1-2\vert 2z\vert-\vert 2z\vert^2)}, & \vert z\vert>1/2.
\end{array} \right.
\end{align*}

Let $\phi_n$ be a sequence of radially symmetric non-negative mollifiers with compact support, i.e. for all $n\geq 0$
\begin{itemize}
\item[(i)] $\phi_n(z)=\rho_n(\vert z\vert)$ for some $\rho_n:[0,+\infty)\to[0,+\infty)$.
\item[(ii)] $\phi_n(z)\geq 0$ for $z$ and $\phi_n\in C^\infty_0(\C)$.
\item[(iii)] $\int_\C\phi_n(z)dA(z)=1$.
\item[(iv)] $\text{supp}(\phi_n)\subset r_n\Di$, $\lim_{n\to \infty}r_n=0$ and $\lim_{n\to \infty}\phi_n=\delta_0$. 
\end{itemize}
Consider the convolution
\begin{align*}
\A_n(z)=\int_{\R^2}\A(z-w)\phi_n(w)dA(w)
\end{align*}
Then $\A_n\in C^\infty(2\overline{\Di})$. Furthermore,

\begin{align*}
&\vert \dv_{\overline{z}}\A_n(z)\vert \leq \int_{2\Di}\vert \A_{\overline{z}}(z-w)\vert \phi_n(w)dA(w)\leq \frac{1}{4}\int_{2\Di}\phi_n(w)dA(w)=\frac{1}{4},\\
&\mathfrak{R}e [\dv_{\overline{z}}\A_n(z)]=  \int_{2\Di}\vert \A_{\overline{z}}(z-w)\vert \phi_n(w)dA(w)\geq \frac{1}{2}\int_{2\Di}\phi_n(w)dA(w)=\frac{1}{2},
\end{align*}
and so $\mathcal{A}_n$ are $\delta$-monotone on $2\Di$.

Finally, for $\vert z\vert\leq 1/4$ 
\begin{align*}
\dv_{\overline{z}}\A_n(z)=\int_{2\Di}\A_w(w)\phi_n(z-w)dA(w)=\int_{r_n\Di+z}\A_w(w)\phi_n(z-w)dA(w)=0
\end{align*}
for $n$ sufficiently large and $\vert z\vert<1/4$. Thus $\nu_n(z)=0$ for $\vert z\vert<1/8$. Moreover, since $\dv_z\A_n(z)$ is real for all $n$, it follows that $\mathfrak{I}m[\dv_z\A_n(z)]=0$ which is equivalent to $\text{curl}\,\A_n(z)=0$. Thus there exists a convex function $\sigma_n: 2\Di\to \R$ such that $\nabla \sigma_n(z)=\A_n(z)$.  In this case one {\bf cannot} use the second equation in \eqref{eq:2} to do the reduction.

\subsection*{Acknowledgements}
Erik Duse was supported by the Knut and Alice Wallenberg Foundation grant KAW 2015.0270. The author thanks Xiao Zhong for reading the manuscript and providing useful suggestions to improve the presentation of the paper.

\bibliographystyle{alpha}

{\sc Erik Duse}, Department of Mathematics and Statistics, KTH,  Stockholm, Sweden \texttt{duse@kth.se}

\end{document}